\documentclass[11pt]{article}
\usepackage{amssymb}
\usepackage{latexsym}
\newcommand{\be}{\begin{equation}}
\newcommand{\ee}{\end{equation}}
\newcommand{\bea}{\begin{eqnarray}}
\newcommand{\eea}{\end{eqnarray}}
\newcommand{\bean}{\begin{eqnarray*}}
\newcommand{\eean}{\end{eqnarray*}}
\newcommand{\benu}{\begin{enumerate}}
\newcommand{\eenu}{\end{enumerate}}
\newcommand{\edo}{\end{document}}

\newcommand{\barr}{\begin{array}}
\newcommand{\earr}{\end{array}}

\newcommand{\CC}{{\mathbb C}}

\newcommand{\HH}{{\mathbb H}}

\newcommand{\RR}{{\mathbb R}}

\newcommand{\Ric}{\mathrm{Ric}}
\newcommand{\mcr}{\mathcal{R}}

\newtheorem{prop}{Proposition}[section]
\newcommand{\ecke}{\;_-\!\rule{0.2mm}{0.2cm}\;\;}

\begin{document}

\title{Conformal Killing spinors and special geometric structures in 
Lorentzian geometry - a survey}

\author{Helga Baum}



\maketitle
\begin{abstract}
This paper is a survey on special geometric structures that admit conformal 
Killing spinors based on lectures, given at the ``Workshop on Special 
Geometric Structures in String Theory'', Bonn,
September 2001 and at ESI, Wien, November 2001. We discuss the case of 
Lorentzian signature and explain which geometries occur up to dimension 6.
\end{abstract}

\tableofcontents

\section{Introduction}

A classical object in differential geometry are conformal Killing fields. These are by definition 
infinitesimal conformal symmetries i.e. the flow of such vector fields 
preserves the conformal class of the metric. The number of linearely independing conformal Killing fields measures the degree of conformal symmetry of the manifold. This number is bounded by $\frac{1}{2}(n+1)(n+2)$, where $n$ is the dimension of the manifold. It is the maximal one if the manifold is conformally flat. S. Tachibana and T. Kashiwada (cf. \cite{Tachibana/Kashiwada:69}, \cite{Kashiwada:68}) introduced a generalisation of conformal Killing fields, the  conformal Killing forms (or twistor forms). Conformal Killing forms are solutions of a conformally invariant twistor type equation on differential forms. They were studied in General Relativity mainly from the local viewpoint in order to integrate the equation of motion (e.g. \cite{Penrose/Walker:70}), furthermore they were used to obtain symmetries of field equations (\cite{Benn/ua1:97}, \cite{Benn/ua2:97}). 
Recently, U. Semmelmann (\cite{Semmelmann:01}) started to discuss global properties of conformal Killing forms in Riemannian geometry. Another generalisation of conformal Killing vectors is that of conformal Killing spinors (or twistor spinors), which are solutions of the conformally invariant twistor equation on spinors introduced 
by R. Penrose in General Relativity (cf. \cite{Penrose/Rindler:86}).  Whereas conformal Killing fields are classical symmetries, conformal Killing spinors define infinitesimal symmetries on supermanifolds (cf. \cite{Alekseevski/Cortes/ua:98}). Special kinds of such spinors, parallel and special Killing spinors, occur in supergravity and string theories. In 1989 A. Lichnerowicz and T. Friedrich started a systematic study of twistor spinors in conformal Riemannian geometry. Whereas  the global structure of Riemannian manifolds admitting conformal Killing spinors is quite well understood (cf. e.g. \cite{Lichnerowicz1:88}, \cite{Lichnerowicz2:88},
\cite{Lichnerowicz1:89}, \cite{Friedrich:89} \cite{Lichnerowicz2:90}, \cite{Baum/Friedrich/ua:91},
\cite{Habermann:90}, \cite{Habermann:93}, \cite{Habermann:94}, \cite{Habermann:96},
\cite{Kuehnel/Rademacher:94},
\cite{Kuehnel/Rademacher1:96}, \cite{Kuehnel/Rademacher1:97},
\cite{Kuehnel/Rademacher:98}), the state of art in its origin, Lorentzian geometry, is far from being satisfactory. J. Lewandowski (\cite{Lewandowski:91}) described the local normal forms of 4-dimensional spacetimes with zero free twistor spinors. His results indicated that there are interesting relations between conformal Killing spinors, different global contact structures and Lorentzian geometry, that should be discovered. In the present survey we discuss some of these structures. Since there is some interest in several kinds of Killing spinors on Lorentzian manifolds by physicists working in string theory, our investigations and methods, motivated mainly from geometry, could be, perhaps, of some use in physics.

\section{Conformally invariant operators on spinors}

In this section we will define the kind of spinor fields we are interested in and discuss two applications 
in Riemannian geometry. \\
Let $(M^n,g)$ be a semi-Riemannian spin manifold of dimension $n \ge 3$. We denote by $S$
 the spinor bundle and by $\mu :T^*M \otimes S \to S$ the 
Clifford multiplication. The 1-forms with values in the spinor bundle 
 decompose into two subbundles
\begin{displaymath}
T^* M \otimes S = V \oplus Tw,
\end{displaymath}
where $V$, being the orthogonal complement to the ``twistor bundle'' $Tw := Ker \mu\;$, is isomorphic to $S$.
 Usually, we identify $TM$ 
and $T^*M$ using the metric $g$. \\
We obtain two differential operators of first order by composing the spinor
derivative $\nabla^S$ with the orthogonal projections onto each of these
subbundles,\\
the {\em Dirac operator} $D$
\begin{displaymath}
D: \Gamma (S) \stackrel{\nabla^S}{\longrightarrow} \Gamma (T^* M \otimes S)=
\Gamma (S \oplus Tw) \stackrel{pr_S}{\longrightarrow}\Gamma (S)
\end{displaymath}
and the {\em twistor operator} $P$
\begin{displaymath}
P: \Gamma (S) \stackrel{\nabla^S}{\longrightarrow} \Gamma (T^*M \otimes S) =
\Gamma (S \oplus Tw) \stackrel{pr_{Tw}}{\longrightarrow} \Gamma (Tw).
\end{displaymath}
Locally, these operators are given by the following formulas
\begin{eqnarray*}
D \varphi &=& \sum\limits^n_{i=1} \sigma^i \cdot \nabla^S_{s_i} \varphi\\
P \varphi &=& \sum\limits^n_{i=1} \sigma^i \otimes (\nabla_{s_i}^S \varphi +
\frac{1}{n} s_i \cdot D \varphi),
\end{eqnarray*}

where  $(s_1 , \ldots , s_n)$ is a local orthonormal basis, $(\sigma^1, 
\ldots, \sigma^n)$ its dual and $\cdot$ denotes the Clifford multiplication.
 Both operators are conformally
covariant. More exactly, if $\tilde{g} = e^{2 \sigma } g$ is a conformal
change of the metric, the Dirac and the twistor operator satisfy
\begin{eqnarray*}
D_{\tilde{g}} &=& e^{- \frac{n+1}{2} \sigma} D_g e^{\frac{n-1}{2} \sigma}\\
P_{\tilde{g}} &=& e^{- \frac{\sigma}{2}} P_g e^{- \frac{\sigma}{2}}.
\end{eqnarray*}
The spinor fields we are interested in are the solutions of the conformally invariant
equations $D \varphi =0$ and $P \varphi =0$. A spinor field $\varphi \in
\Gamma (S)$ is called
\begin{tabbing}
 xx \= twistor spinor (conformal Killing spinor)  \= x$\Longleftrightarrow$ xxx\= twistor spinor (conformal Killing spinor
  \= \kill
\> {\em harmonic spinor} \> $\Longleftrightarrow$ \> $D\varphi = 0$\\[0.1cm]
\> {\em twistor spinor (conformal Killing spinor)} \> $\Longleftrightarrow$  \> $P \varphi =0$ \\[0.1cm]
\> {\em parallel spinor} \> $\Longleftrightarrow$ \> $ \nabla^S \varphi =0$\\[0.1cm]
\> {\em geometric Killing spinor} \> $\Longleftrightarrow $ \>  $\nabla_X^S \varphi = \lambda X \cdot \varphi ,\; 
\lambda \in \mathbb{C}$\\
\> {\em with Killing number} $\lambda$ \> \> 

\end{tabbing}

The local formula for the twistor operator shows another characterisation of a twistor spinor: 
A spinor field $\varphi \in \Gamma(S)$ is a twistor spinor if and only if 
\begin{displaymath}
\fbox{$\nabla^S_X \varphi + \frac{1}{n} X \cdot D\,\varphi = 0$. }
\end{displaymath}
Obviously, parallel and Killing spinors are special kinds of twistor spinors.
Each twistor spinor that satisfies the Dirac equation $D \varphi = \mu
\varphi$ is a Killing spinor. On Riemannian manifolds each
twistor spinor without zeros can be transformed by a conformal change of the
metric into a Killing spinor. 
\\
Twistor spinors appeared naturally in Riemannian geometry. 
For a motivation, 
let me explain two well known examples, where it was usefull to know 
 structure results for manifolds with special kinds of twistor spinors.\\[0.2cm]
{\bf 1. Eigenvalue estimates for the Dirac operator on closed Riemannian
manifolds with positive scalar curvature $R$ (limiting case)}\\[0.1cm]
The Dirac operator $D$ on a closed Riemannian spin manifold is elliptic and
essentially selfadjoint. Hence, its spectrum contains only real eigenvalues
of finite multiplicity. One is interested in lower estimates for the smallest
eigenvalue. It is easy to prove that the Dirac operator $D$ and the twistor
operator $P$ are related by 
\begin{displaymath}
D^2 = \frac{n}{n-1} (P^* P + \frac{1}{4} R) \ . 
\end{displaymath}
If $\lambda$ is an eigenvalue of $D$ and $D \varphi = \lambda \varphi$, 
integration of
\begin{displaymath}
\langle D^2 \varphi, \varphi \rangle = \frac{n}{n-1} \langle P^* P \varphi , \varphi \rangle +
\frac{1}{4}\langle R \varphi, \varphi \rangle
\end{displaymath}
yields the estimate (\cite{Friedrich:80}) 
\begin{equation} \label{gl01}
\lambda^2 \ge \frac{1}{4} \frac{n}{n-1} \cdot \min_{x \in M} R (x),
\end{equation}
where equality holds if and only if the eigenspinor $\varphi$ is a twistor spinor. 
Hence, the discussion of the limiting case (i.e. equality holds in
(\ref{gl01})) leads to the problem to find all Riemannian structures which 
admit real Killing spinors (i.e. Killing spinors with real Killing number).
 This problem was studied by several authors 
\footnote{Friedrich, Grunewald, Kath, Hijazi, Baum, B\"ar, Duff, Nillson, 
Pope, Nieuwenhuizen and many others)} between 1980 and 1993 and is now 
essentially solved. The key observation in the understanding of geometric Killing spinors was made
by Ch. B\"ar in 1993 (\cite{Baer:93}). He remarked, that a Riemannian manifold
$(M^n,g)$ has real Killing spinors 
if and only if the metric cone
$C_+ (M) = ( \RR^+ \times M, dt^2 + t^2 g )$
has parallel spinors. The cone is irreducible in case $M^n$ is not the standard sphere.
Then, to the cone, one can apply a result 
of McK. Wang (\cite{Wang:89}) that describes the possible reduced holonomy groups
$\mathrm{Hol}_0$ of irreducible manifolds with parallel spinors and is able to derive 
the corresponding geometries of $M$. 
The result is given in the following table. A compact simply connected Riemannian manifold
$(M^n,g)$ admits real Killing spinors if and only if
\begin{center}
\begin{tabular}{|c|c|c|c|c|c|}
\hline &&&&& \\[-0.2cm]
$\mathrm{Hol}(C_+(M))$ & $\mathrm{SU} (\frac{n}{2})$ & $\mathrm{Sp} (
\frac{n}{4})$ & $G_2$ & $\mathrm{Spin} 7$ & $1$ \\ &&&&& \\[-0.2cm] \hline &&&&&\\[-0.2cm]
$M^n$ & Einstein-Sasaki & 3-Sasaki & nearly K\"ahler & nearly parallel & $S^n$\\
&&& non K\"ahler & $G_2$-structure & \\[-0.2cm] &&&&& \\ \hline
\end{tabular}
\end{center}
\ \\
{\bf Rigidity Theorems}\\[0.2cm]
As a second example I want to explain the Rigidity Theorem for asymptotically locally hyperbolic (ALH) manifolds,
proved by Andersson and Dahl (\cite{Andersson/Dahl:98}): An ALH manifold
is a complete Riemannian manifold with an asymptotically locally hyperbolic
end $\mathbb{H}^n /\Gamma$. Andersson and Dahl proved that each ALH
manifold of dimension $n$ with scalar curvature $R$ bounded by
$R \ge - n(n-1)$ is isometric to the hyperbolic space $\mathbb{H}^n$. 
The hyperbolic space $\mathbb{H}^n$ admits imaginary Killing spinors (i.e. the Killing number is purely imaginary).
Using such Killing spinor on the asymptotic end of the ALH manifold one can
construct an imaginary Killing spinor on $M^n$ itself. A complete Riemannian
spin manifold $(M^n,g)$ admits imaginary Killing spinors if and only if 
$(M^n,g)$ is isometric to a warped product
$(\RR \times F , dt^2 + e^{t\mu} g_F )$ , 
where $(F,g_F)$ is a complete Riemannian manifold with parallel spinors 
(\cite{Baum2:89}). But, being asymptotically locally hyperbolic, such a warped product has to be
 the hyperbolic space $\mathbb{H}^n$.\\[0.4cm]
Whereas twistor spinors on Riemannian manifolds are rather well studied
 \footnote{e.g., by Lichnerowicz, Friedrich, K.Habermann, K\"uehnel, Rademacher},
the state of art in Lorentzian geometry, where the twistor equation originally came from, 
is far from being satisfactory. For that reason, 
in this survey we want to consider twistor spinors in Lorentzian geometry 
and discuss the following problems:
\begin{enumerate}
\item {\em Which Lorentzian geometries admit twistor spinors?}
\item {\em How the properties of twistor spinors are related to the geometric 
structures where they can occur?}
\end{enumerate}

\section{Basic facts on twistor spinors}

There is a fundamental difference between the two conformally invariant
equations $D \varphi =0$ and $P \varphi =0$. Whereas the dimension of 
space of harmonic spinors can become arbitrary large, the dimension of the space of twistor spinors
is always bounded by $2\,\mathrm{rank} S$. It holds even more: the twistor equation $P \varphi =0$ can 
be viewed as a parallelity equation with respect to a suitable covariant 
derivative in a suitable bundle. Let us explain this fact (well known for a 
long time, cf. \cite{Penrose/Rindler:86}, \cite{Baum/Friedrich/ua:91})
more in detail.
Let $\mathcal{R}$ and $W$ denote the curvature tensor and the Weyl tensor of $(M^n,g)$, respectively, considered
 as selfadjoint maps on the space
of 2-forms
$\mathcal{R}, W : \Lambda^2 M \to \Lambda^2 M\,$.
$\,\mathrm{Ric}$ denotes the Ricci tensor of $(M^n,g)$ considered here as
$(1,1)$-tensor or as $(2,0)$-tensor as it is needed.
In conformal geometry there are 2 further curvature tensors that play an 
important role, the Rho tensor $K$
\begin{displaymath}
K(X) := \frac{1}{n-2} \Big( \frac{R}{2(n-1)} X - \Ric (X) \Big) \quad,
X \in TM
\end{displaymath}
and the Cotton-York tensor
\begin{displaymath}
C(X,Y):= (\nabla_X K)(Y) - (\nabla_Y K)(X) , \quad X,Y \in TM.
\end{displaymath}
Let us consider the double spinor bundle $E=S \oplus S$ with the 
following covariant derivative
\begin{displaymath}
\nabla_X^E := \left( \begin{array}{cc}
\nabla_X^S & \frac{1}{n} X \cdot \\
- \frac{n}{2} K(X) \cdot & \nabla^S_X \end{array} \right) 
\end{displaymath}
The curvature of this derivative is given by
\begin{equation}  \label{gl02}
\mcr^{\nabla^E} (X,Y) \left( \begin{array}{c} \varphi \\ \psi \end{array}
\right) = \frac{1}{2} \left( \begin{array}{l}
W(X \wedge Y) \cdot \varphi\\
W(X \wedge Y) \cdot \psi - n C(X,Y) \cdot \varphi
\end{array}\right) 
\end{equation}
(cf. \cite{Baum/Friedrich/ua:91}, chap. 1).
Using integrability conditions for twistor spinors one obtains the following
correspondence between twistor spinors and parallel spections in $(E, 
\nabla^E)$:

\begin{prop} \label{prop-1}
Let $(M^n,g)$ be a semi-Riemannian spin manifold and $\varphi \in \Gamma
(S)$. Then\\[-0.7cm]
\begin{enumerate}
\item 
If $\varphi$ is a twistor spinor, then $\nabla^E \left( \begin{array}{c} \varphi \\
 D\varphi \end{array} \right)=0. $ \\[-0.6cm]
\item 
If $\nabla^E \left( \begin{array}{c} \varphi \\ \psi \end{array} \right) 
=0$, then $\varphi$ is a twistor spinor and $\psi = D \varphi$.
\end{enumerate}
\end{prop}

A semi-Riemannian manifold is (locally) conformally flat if and only if $W=0$ (hence $C=0$)
 in dimension $n \geq 4$ or $C=0$ in
dimension $n=3$, where $W=0$ is always satisfied. Therefore, by formula
(\ref{gl02}) the spin manifold $(M^n,g)$ is conformally flat if and only if
the curvature $\mcr^{\nabla^E}$ of $\nabla^E$ vanishes. Together with Proposition \ref{prop-1} one 
obtains the maximal possible number of linearely independent twistor spinors, which is attained on conformally flat manifolds as in the case of conformal vector fields.

\begin{prop}  \label{prop-2} 
\begin{enumerate}
\item The dimension of the space of twistor spinors is a conformal 
invariant and bounded by
\begin{displaymath}
\dim \ker P \le 2 \,\mathrm{rank} S = 2^{[\frac{n}{2}]+1} =: d_n.
\end{displaymath}
\item If  $\dim \ker P =d_n$, then
$(M^n,g)$ is conformally flat.\\[-0.7cm]
\item If $(M^n,g)$ is simply connected and conformally flat, then $\dim \ker 
P= d_n$.
\end{enumerate}
\end{prop}

For example, all simply connected space forms $\RR^n_k, \tilde{\HH}^n_k
, \tilde{S}^n_k$ admit the maximal number of linearly independent twistor
spinors.

\section{Twistor spinors on Lorentzian spin manifolds}

Now, we restrict our attention to the case of Lorentzian signature $(-,+ 
\ldots +)$. Let $(M^n,g)$ be an oriented and time-oriented Lorentzian spin 
manifold. On the spinor bundle $S$ there exists an indefinite non-degenerate
inner product $\langle \cdot ,\cdot \rangle$ such that
\begin{eqnarray*}
\langle X,  \cdot \varphi  , \psi \rangle & =& \langle \varphi , X \cdot \psi
\rangle \quad \quad \mbox{and}\\
X(\langle \varphi , \psi \rangle )&= &\langle \nabla^S_X \varphi , \psi 
\rangle + \langle \varphi , \nabla^S_X \psi \rangle,
\end{eqnarray*}
for all vector fields $X$ and all spinor fields $\varphi, \psi$. Each
spinor field $\varphi \in \Gamma (S)$ defines a vector field $V_{\varphi}$
on $M$, the so-called Dirac current, by
\begin{equation} \label{gl03}
g (V_{\varphi} , X) := - \langle X \cdot \varphi , \varphi \rangle . 
\end{equation}
A motivation to call twistor spinors conformal Killing spinors is the 

\begin{prop} \label{prop-3}
Let $\varphi$ be a spinor field on a Lorentzian spin manifold $(M^n,g)$
with Dirac current $V_{\varphi}$. Then\\[-0.7cm]
\begin{enumerate}
\item $V_{\varphi}$ is causal and future-directed.\\[-0.7cm]
\item The zero sets of $\varphi$ and $V_{\varphi}$ coincide.\\[-0.7cm]
\item If $\varphi$ is a twistor spinor, $V_{\varphi}$ is a conformal Killing
field.
\end{enumerate}
\end{prop}
Now, let us discuss 3 types of special Lorentzian geometries that admit conformal 
Killing spinors.

\subsection{Brinkman spaces with parallel spinors}

A Lorentzian manifold is called \emph{Brinkman space} if it admits a 
non-trivial  lightlike parallel vector field. Let us consider two examples of such spaces.\\

{\bf Example 1}: pp-manifolds.\\[0.1cm]
A \emph{pp-manifold} is a Brinkman space with the additional curvature 
property
\begin{displaymath}
\mbox{Trace}\,_{(3,5),(4,6)} \mcr \otimes \mcr =0 . 
\end{displaymath}
Equivalently, pp-manifolds can be characterised as those Lorentzian
manifolds $(M^n,g)$, where the metric has the following local normal form
depending only on one function $f$ of $(n-1)$ variables
\begin{displaymath}
g= dt\, ds +f(s, x_1 , \ldots , x_{n-2}) ds^2 + \sum\limits^{n-2}_{i=1}
dx_i^2 . 
\end{displaymath}
(cf. \cite{Schimming:74}).
In terms of holonomy, pp-manifolds can be characterised as those 
Lorentzian manifolds for which the  holonomy algebra
is an abelian ideal of $\mathfrak{so} (1, n-1)$. (see \cite{Leistner:01}). 
Using the latter fact one can easily prove that for each simply connected
pp-manifold
\begin{displaymath}
\dim \ker P \ge \frac{d_n}{4} .
\end{displaymath}

Furthermore, on pp-manifolds each twistor spinor is parallel. An 
important examples of geodesically complete pp-manifolds are the Lorentzian
symmetric spaces with solvable transvection group which occur also as a
special model for a certain string equation in $10+1$-dimension (cf. the 
talk of J. Figueroa O`Farrill at this Workshop or \cite{Figueroa1:01}, \cite{Figueroa2:01}).\\

{\bf Example 2}: Brinkman spaces with special K\"ahler flag.\\
Let $(M^n,g)$ be a Brinkman space with the  lightlike parallel vector field
$V$. $V$ defines a flag $\,\RR V \subset V^{\perp} \subset TM \,$ in $TM$,
where $V^{\perp} =\{ Y \in TM \,|\, g(V,Y)=0 \}$. We equip the bundle $E:= V^{\perp} /\RR V$  
with the 
positive definite inner product $\tilde{g}$ induced by $g$ and the  metric
connection $\tilde{\nabla}$ induced by the Levi-Civita connection of $g$. We 
call the flag $\RR V \subset V^{\perp} \subset TM$ a \emph{special
K\"ahler flag}, if the bundle $E$ admits an orthogonal almost complex structure
$J:E \to E,\, J^2 = - id\, $, such that
$\;\mbox{Trace}\; (J \circ \mcr^{\tilde{\nabla}} (X,Y)) =0\,$
for all $X,Y \in TM$. 
It was proved by I. Kath in \cite{Kath:99} that $(M^{2m},g)$ is a 
Brinkman space with special K\"ahler flag if and only if $(M^{2m},g)$ has
pure parallel spinors.

\subsection{Twistor spinors on Lorentzian Einstein Sasaki manifolds}

An odd-dimensional Lorentzian manifold $(M^{2m+1}, g ; \xi)$ equipped
with a vector field $\xi$ is called \emph{Lorentzian Sasaki
manifold} if\\[-0.7cm]
\begin{enumerate}
\item $\xi$ is a timelike Killing field with $g( \xi, \xi )=-1.$\\[-0.7cm]
\item The map $J:= - \nabla \xi :TM \to TM$ satisfies \\[-0.7cm]
\begin{eqnarray*}
J^2 X &=& -X -g (X, \xi ) \xi \quad \quad \mbox{and}\\
(\nabla_X J)(Y) &=& - g(X,Y) \xi + g (Y, \xi )X.
\end{eqnarray*}
\end{enumerate}

Let us consider the metric cone $\,C_- (M) := (\RR^+ \times M , -dt^2 + t^2 g)\,$ with timelike cone axis over 
$(M,g)$. The cone metric has signature $(2,2m)$. Then the following relations
between properties of $M$ and those of its cone are easy to verify 

\begin{tabbing}
xxxxxxxxxxxxxxxxxxxxxxxxxxxxxxxx\= xxxxxxxx\= \kill
$(M^{2m+1},g;\xi)$ \> \>  cone $C_-(M)$ \\[0.2cm]  
 Lorentzian Sasaki \> $\Longleftrightarrow$ \> (pseudo) K\"ahler\\ 
 Lorentzian Einstein Sasaki ($R<0$) \> $\Longleftrightarrow$ \>  
Ricci-flat and (pseudo) K\"ahler\\
 Lorentzian Einstein Sasaki ($R<0$) \> $\Longleftrightarrow$\>  $\mathrm{Hol}_0 (C_- (M)) \subset
SU(1,m)$
\end{tabbing}

The standard example for regular Lorentzian Einstein Sasaki manifolds are $S^1$- 
bundles over Riemannian K\"ahler Einstein spaces of negative scalar curvature:
Let $(X^{2m},h)$ be a Riemannian K\"ahler Einstein spin manifold of scalar
curvature $R_X <0$ and let $(M^{2m+1}, \pi , X ; S^1)$ denote the 
$S^1$-principal bundle associated to the square root $\sqrt{\Lambda^{m,0} X}$
of the canonical bundle of $X$ given by the spin structure. Furthermore,
let $A$ be the connection on $M$ induced by the Levi-Civita
connection of $(X,h)$. Then
\begin{displaymath}
g := \pi^* h - \frac{16m}{(m+1) R_X} A \circ A
\end{displaymath}

defines a Lorentzian Einstein Sasaki metric on the spin manifold 
$M^{2m+1}$. Lorent\-zian Einstein Sasaki manifolds admit a special kind
of twistor spinors.

\begin{prop} (\cite{Kath:99}, \cite{Bohle:00}) \label{prop-4}
Let $(M,g)$ be a simply connected Lorentzian Einstein Sasaki manifold. Then
$M$ is spin and admits a twistor spinor $\varphi$ on $M$ such that\\[-0.7cm]
\begin{itemize}
\item[a)] $V_{\varphi}$ is a timelike Killing field with
$g(V_{\varphi} , V_{\varphi})=-1$\\[-0.7cm]
\item[b)] $V_{\varphi}  \cdot \varphi = - \varphi$\\[-0.7cm]
\item[c)] $\nabla_{V_{\varphi}} \varphi = - \frac{1}{2} i \varphi$.
\end{itemize}
Conversally, if $(M,g)$ is a Lorentzian spin manifold with a twistor spinor 
satisfying a), b) and c). Then $\xi := V_{\varphi}$ is a Lorentzian
Einstein Sasaki structure on $(M,g)$.
\end{prop}

\subsection{Twistor spinors on Fefferman spaces}

Fefferman spaces are Lorentzian manifolds which appear in the frame work
of CR geometry. Let us first explain the necessary notations from CR 
geometry. Let $N^{2m+1}$ be a smooth oriented manifold of odd dimension
$2m+1$. A \emph{CR structure on $N$} is a pair $(H,J)$ where\\[0.2cm]
\hspace*{0.1cm} 1.  $H \subset TM$ is a real $2m$-dimensional subbundle.\\[0.1cm]
\hspace*{0.1cm} 2.  $J:H \to H$ is an almost complex structure on $H$, $J^2 =-Id$.\\[0.1cm]
\hspace*{0.1cm} 3.  If $X,Y \in \Gamma (H)$, then $[JX,Y]+[X,JY] \in \Gamma (H)$ and\\[0.1cm]
    \hspace*{1cm}  $J([JX,Y]+[X,JY]) - [JX, JY]+[X,Y] \equiv 0$ (integrability condition).\\

Standard examples of CR manifolds are the following:\\[-0.5cm]
\begin{itemize}
\item Real hypersurfaces $N$ of a complex manifold $(M,J_M)$. The CR structure is given by 
$H:=TN \cap J_M (TN) , J:=J_M \big|_{{\Large H}}$.\\[-0.5cm]
\item Riemannian Sasaki manifolds $(N,g, \xi)$. The CR structure is given by $H:= \xi^{\perp}$
and $J:= - \nabla \xi$.\\[-0.6cm]
\item Heisenberg manifolds $\mathcal{H}e(m) =He(m)/\Gamma$, where
$He (m)$ is the Heisenberg group of matrices
\begin {displaymath}
He(m) =\left\{ \left( \begin{array}{ccc}
1& x^t &z\\
0&I_m &y\\
0&0&1
\end{array} \right) \,\Bigg|\, x,y \in \RR^m , z \in \RR \right\}
\end{displaymath}
and $\Gamma$  is a discret lattice in $He (m)$.
\end{itemize}

Let $(N,H,J)$ be a CR manifold. In order to define Fefferman spaces we fix a contact form 
$\theta \in \Omega^1 (N)$ on $N$ such that
$\theta|_H = 0$. Let us denote by $T$ the Reeb vector field of $\theta$
which is defined by the conditions $\theta (T) =1$ and  $T \ecke d \theta =0$. In the 
following we suppose that the Levi form $L_{\theta} :H \times H \to \RR$
\begin{displaymath}
L_{\theta} (X,Y) := d \theta (X,JY)
\end{displaymath}
is positive definite. Then $(N,H, J , \Theta)$ is called a 
{\em strictly pseudoconvex manifold}. The tensor field $g_{\theta} := L_{\theta}
+ \theta \circ \theta$ defines a Riemannian metric on $N$. There is a special
 metric covariant derivative on a strictly pseudoconvex manifold, the 
{\em Tanaka-Webster connection} $\nabla^W : \Gamma (TN) \to \Gamma (TN^* 
\otimes TN)$ uniquely defined by the conditions
\begin{eqnarray*}
\nabla^W g_{\theta} &=& 0\\
Tor^W (X,Y) &=& L_{\theta} (JX,Y) \cdot T\\
Tor^W (T,X) &=& - \frac{1}{2} ([T,X]+J[T,JX])
\end{eqnarray*}
for $X,Y \in \Gamma (H)$. This connection satisfies $\nabla^W J= 0$ and
$\nabla^W T=0$ (cf. \cite{Tanaka:75}, \cite{Webster:78}). Let us denote by
$T_{10} \subset TN^{\CC}$ the eigenspace of the complex extension of $J$ on $H^{\CC}$
to the eigenvalue $i$. Then $L_{\theta}$ extends to a hermitian form on $T_{10}$ by 
\[L_{\theta} (U,V) := - i\,d\theta (U, \bar{V}), \;\;\;U,V \in T_{10}.\]
For a complex 2-form $\omega \in \Lambda^2 N^{\mathbb{C}}$ we denote by
$\mbox{trace}_{\theta} \omega$ the $\theta$-trace of $\omega$:
\begin{displaymath}
\mbox{trace}_{\theta} \omega := \sum\limits^m_{\alpha =1} \omega (Z_{\alpha} , 
\bar{Z}_{\alpha} ), 
\end{displaymath}
where $(Z_1 , \ldots , Z_m)$ is an unitary basis of $(T_{10}, L_{\theta})$.
Let ${\cal R}^W$ be the $(4,0)$-curvature tensor of the Tanaka-Webster
connection $\nabla^W$ on the compexified tangent bundle of $N$
\begin{displaymath}
{\cal R}^W (X,Y,Z,V) := g_{\theta} (([\nabla^W_X , \nabla^W_Y]-
\nabla^W_{[X,Y]}) Z, \bar{V}) . 
\end{displaymath}
and let us denote by
\begin{displaymath}
\Ric^W := \mbox{trace}^{(3,4)}_{\theta} := \sum\limits^m_{\alpha =1} {\cal R}^W
(\cdot , \cdot , Z_{\alpha} , \bar{Z}_{\alpha})
\end{displaymath}
the \emph{Tanaka-Webster Ricci curvature} and by $R^W := \mbox{trace}_{\theta}
\Ric^W$ the \emph{Tanaka-Webster scalar curvature}. The Ricci curvature $\Ric^W$ is a 
$(1,1)$-form on $N$ with $\Ric^W (X,Y) \in i \RR$ for real vectors $X,Y \in 
TN$. The scalar curvature $R^W$ is a real function.\\
Now, let us suppose that $(N^{2m+1}, H, J, \theta)$ is a strictly
pseudoconvex spin manifold. The spin structure of $(N, g_{\theta})$ defines 
a square root $\sqrt{\Lambda^{m+1,0} N}$ of the canonical line bundle
\begin{displaymath}
\Lambda^{m+1,0} N:= \{ \omega \in \Lambda^{m+1} N^C \,|\, V \ecke 
\omega =0 \,\,\,\forall V \in \bar{T}_{10} \} . 
\end{displaymath}

We denote by $(F, \pi , N)$ the $S^1$-principal bundle associated to 
$\sqrt{\Lambda^{m+1,0} N}$. 
If one fixes a connection form $A$ on $F$ and the corresponding decomposition
of the tangent bundle $TF = ThF \oplus TvF= H^* \oplus \RR T^* \oplus TvF$
into the horizontal and vertical part, then a Lorentzian metric $h$ is
defined by
\begin{displaymath}
h := \pi^* L_{\theta} - i\,c \,\pi^*\theta \circ A , 
\end{displaymath}
where $c$ is a non-zero real number.
The Fefferman metric is a metric of the latter type, where the  choice of $A$ and $c$ is done in
such a way that the conformal class  of $h$ does not depend on the
pseudohermitian form $\theta$. Such a choice can be made with the
connection
\begin{displaymath}
A_{\theta} := A^W - \frac{i}{4(m+1)} R^W \cdot \theta , 
\end{displaymath}
where $A^W$ is the connection form of $F$ defined by the Tanaka-Webster
connection $\nabla^W$. The curvature form of $A^W$ is $\Omega^{A^W} = - 
\frac{1}{2} \Ric^W$. Then
\begin{displaymath}
h_{\theta} := \pi^* L_{\theta} - i \,\frac{8}{m+2} \pi^* \theta \circ
A_{\theta}
\end{displaymath}
is a Lorentzian metric such that the conformal class $[h_{\theta}]$ is an
invariant of the CR-structure $(N,H,J)$. The metric $h_{\theta}$ is $S^1$-invariant, the fibres 
of the $S^1$-bundle are lightlike. We call
$(F^{2m+2}, h_{\theta})$ with its canonically induced spin structure
\emph{Fefferman space of the strictly pseudoconvex spin manifold $(N,H,J,
\theta)$}.\\[0.1cm]
The Fefferman metric was first discovered by C. Fefferman for the case
of strictly pseudoconvex hypersurfaces $N \subset \mathbb{C}^{m+1}$
(\cite{Fefferman:76}), who showed that $N \times S^1$ carries a Lorentzian metric 
whose conformal class is induced by biholomorphisms. The considerations
of Fefferman were extended by Burns, Diederich and Snider (\cite{Burns/Diederich/ua:77}) and
by Lee (\cite{Lee:86}) to the case of abstract (not nessarily embedded)
CR-manifolds. A geometric characterization of Fefferman metrics was given by
Sparling (cf. \cite{Sparling:85}, \cite{Graham:87}).\\
Fefferman spaces admit a special kind of twistor spinors.

\begin{prop} (\cite{Baum:99}) \label{prop-5}
Let $(N,H,J, \theta)$ be a strictly pseudoconvex spin manifold with the
Fefferman space $(F, h_{\theta})$. Then there exist two linearly independent
twistor spinors $\varphi$ on $(F, h_{\theta})$ such that \\[-0.6cm]
\begin{itemize}
\item[a)] $V_{\varphi}$ is a regular lightlike Killing field\\[-0.7cm]
\item[b)] $V_{\varphi} \cdot \varphi =0$\\[-0.7cm]
\item[c)] $\nabla_{V_{\varphi}} \varphi = i\, c\, \varphi$ where
$c \in \RR \backslash \{ 0 \}$.
\end{itemize}
Conversally, if $(M,g)$ is an even dimensional Lorentzian spin 
manifold with a twistor spinor satisfying a), b) and c), then there
exists a strictly pseudoconvex spin manifold $(N,H,J, \theta)$ such 
that its Fefferman space is locally isometric to $(M,g)$. 
\end{prop}

\subsection{Twistor spinors inducing lightlike Killing fields}

As we noticed in Proposition \ref{prop-3} each twistor spinor $\varphi$ 
induces a causal conformal Killing field $V_{\varphi}$. Now we study
the situation where $V_{\varphi}$ is lightlike and Killing. The following result explaines the 
role Fefferman spaces are playing in the set of all special geometries that admit twistor spinors.

\begin{prop} (\cite{Baum/Leitner:02}) \label{prop-6}
Let $(M,g)$ be a Lorentzian spin manifold admitting a twistor spinor 
$\varphi$ such that $V_{\varphi}$ is lightlike and Killing. Then the function
$\Ric (V_{\varphi} , V_{\varphi})$ is constant and non-negative on $M$. 
Furthermore,
\begin{enumerate}
\item $\Ric (V_{\varphi}, V_{\varphi}) >0$ if and only if $(M,g)$ is
locally isometric to a Fefferman space.
\item $\Ric (V_{\varphi} , V_{\varphi}) =0$ if and only if $(M,g)$ is
locally conformal equivalent to a Brinkman space with parallel spinors.
\end{enumerate}
\end{prop}
 
\subsection{Twistor spinors in dimension $n \le 6$}

Finally, we describe all geometric structures of Lorentzian spin manifolds
with twistor spinors that appear locally up to dimension 6.

\begin{prop} (\cite{Lewandowski:91}, \cite{Leitner:01}, \cite{Baum/Leitner:02}) \label{prop-7}
Let $(M^n,g)$ be a Lorentzian manifold with twistor spinors without zeros. 
Then $(M^n,g)$ is locally conformal equivalent to one of the following
kinds of Lorentzian structures.

\begin{center}
\begin{tabular}{|l|l|}
\hline & \\ [-0.2cm]
$n=3$ & pp-manifold \\[-0.2cm] & \\ \hline
&\\[-0.2cm]
$n=4$ & pp-manifold\\
& Fefferman space \\[-0.2cm] & \\ \hline
& \\[-0.2cm]
$n=5$ & pp-manifold\\
& Lorentzian Einstein Sasaki manifold\\
& $(\RR  , - dt^2)\times (N^4,h)$, where $(N^4,h)$ is a Riemannian\\ 
& \hspace{4cm}Ricci-flat
K\"ahler manifold\\[-0.2cm] & \\ \hline
&\\[-0.2cm]
$n=6$ & pp-manifold\\
& Fefferman space\\
& $\RR^{1,1} \times (N^4,h)$, where $(N^4,h)$ is a Ricci-flat K\"ahler 
manifold\\
& Brinkman space with special K\"ahler flag \\[-0.2cm] &\\ \hline
\end{tabular}
\end{center}
\end{prop}

\ \\
\small{

}
\vspace{0.8cm}
{\footnotesize Helga Baum, 
Institut f\"ur  Mathematik,
Humboldt-Universit\"at zu Berlin,
Sitz: Rudower Chaussee 25,
12489 Berlin,
Germany\\
email: baum@mathematik.hu-berlin.de\\}

\end{document}